\documentclass[11pt]{article}

\usepackage{amsmath,amssymb}
\usepackage[mathscr]{eucal}
\usepackage[all]{xy}

\newtheorem{Theorem}{Theorem}[section]
\newtheorem{Lemma}[Theorem]{Lemma}
\newtheorem{Proposition}[Theorem]{Proposition}
\newtheorem{Corollary}[Theorem]{Corollary}

\newtheorem{Definition}[Theorem]{Definition}
\newtheorem{Example}[Theorem]{Example}

\makeatletter \@addtoreset{equation}{section} \makeatother

\title{\bf\LARGE A QUANTUM ANALOGUE OF THE BERNSTEIN FUNCTOR}
\author{S. Sinel'shchikov$^\dagger$\and A. Stolin$^\ddagger$\and L.
Vaksman$^\dagger$}
\date{\sl$^\dagger$Mathematics Division, Institute for Low
Temperature Physics \& Engineering \linebreak 47 Lenin Ave, 61103 Kharkov,
Ukraine \linebreak \linebreak $^\ddagger$Chalmers Tekniska H\"ogskola,
Mathematik, 412 96, G\"oteborg, Sweden}

\begin{document}
\large
\maketitle

\begin{abstract}
We consider Knapp-Vogan Hecke algebras in the quantum group setting. This
allows us to produce a quantum analogue of the Bernstein functor as a first
step towards the cohomological induction for quantum groups.
\end{abstract}

The cohomological induction is one of the most important tools for
producing unitarizable Harish-Chandra modules \cite{KV}. The present paper
demonstrates that this method can be applied to modules over quantum
universal enveloping algebras \cite{Jant} as well. Our way to this circle
of problems was as follows.

The papers \cite{SV, SSV1} introduce quantum analogues of the bounded
symmetric domains and the associated Harish-Chandra modules. One special
case considered in \cite{SSSV2} leads to two geometric realizations for the
ladder representation of the quantum universal enveloping algebra
$U_q\mathfrak{su}_{2,2}$ and to the quantum Penrose transform
\cite{BaEast}. One of the main tools used in \cite{SSSV2} was studying
q-analogues of the \v{C}ech cohomology. Unfortunately this method can be
hardly generalized because it relates to non-trivial problems in
non-commutative algebraic geometry. We hope to overcome these obstacles via
replacing the q-analogues of the \v{C}ech cohomology by q-analogues of the
Dolbeault cohomology. It is well-known that in the classical situation the
Dolbeault cohomology can be constructed algebraically using the so-called
cohomological induction (see \cite{KV}).

It is convenient to work with the Bernstein functor (the projective
Zuckerman functor) for producing unitarizable Harish-Chandra modules. The
definition of the Bernstein functor uses essentially the algebra of
distributions on a real reductive group $G$ with support in a maximal
compact subgroup $K$. The principal obstacle to cope with in this work is
to construct a q-analogue for this algebra.

In the present paper we construct a q-analogue of the algebra of
distributions as above in certain cases. This allows us to construct a
q-analogue of the Bernstein functor.

The authors are grateful to V. Kruglov for assistance in preparation of a
preliminary version of this work and for helpful discussions.

\section{The category \boldmath $C(\mathfrak{l},\mathfrak{l})_q$ and the
algebra $R(\mathfrak{l})_q$}

We will use $\mathbb{C}$ as the ground field.

Let $\mathbf{a}=(a_{ij})_{i,j=1,2,\ldots,l}$, be the Cartan matrix for a
simple complex Lie algebra of rank $l$. Its universal enveloping algebra is
a unital algebra $U$ determined by the generators
$\{H_i,E_i,F_i\}_{i,j=1,2,\ldots,l}$ and the relations \cite[p. 51]{Jant}:
$$H_iH_j-H_jH_i=0,\qquad E_iF_j-F_jE_i=\delta_{ij}H_i,$$
$$
H_iE_j-E_jH_i=a_{ij}E_j,\qquad H_iF_j-F_jH_i=-a_{ij}F_j,\qquad
i,j=1,2,\ldots,l,
$$
and
$$
\sum_{m=0}^{1-a_{ij}}(-1)^m\binom{1-a_{ij}}mE_i^{1-a_{ij}-m}\cdot E_j\cdot
E_i^m=0,
$$
$$
\sum_{m=0}^{1-a_{ij}}(-1)^m\binom{1-a_{ij}}mF_i^{1-a_{ij}-m}\cdot F_j\cdot
F_i^m=0,
$$
for all $i\ne j$.

The Cartan subalgebra spanned by $H_1,H_2,\ldots,H_l$ is denoted by
$\mathfrak{h}$, and $\alpha_1,\alpha_2,\ldots,\alpha_l$,
$$\alpha_j(H_i)=a_{ij},\qquad i,j=1,2,\ldots,l,$$
are simple roots.

There exists a unique collection of coprime positive integers
$d_1,d_2,\ldots,d_l$ such that
$$d_ia_{ij}=d_ja_{ji},\qquad i,j=1,2,\ldots,l.$$
The bilinear form in $\mathfrak{h}^*$ given by
$$(\alpha_i,\alpha_j)=d_ia_{ij},\qquad i,j=1,2,\ldots,l,$$
is positive definite.

Recall the definition of the quantum universal enveloping algebra $U_q$
introduced by V. Drinfeld and M. Jimbo. We assume that
$q\in(0,1)$.\footnote{The purpose of so strong assumptions on $q$ will
become evident while investigating the unitarizability of the
Harish-Chandra modules in question over $U_q$.}

The unital algebra $U_q$ is determined by its generators $K_i$, $K_i^{-1}$,
$E_i$, $F_i$, $i=1,2,\ldots,l$, and the relations
$$K_iK_j=K_jK_i,\qquad K_iK_i^{-1}=K_i^{-1}K_i=1,$$
$$
K_iE_j=q_i^{a_{ij}}E_jK_i,\quad K_iF_j=q_i^{-a_{ij}}F_jK_i,\quad
E_iF_j-F_jE_i=\delta_{ij}(K_i-K_i^{-1})/(q_i-q_i^{-1}),
$$
$$
\sum_{s=0}^{1-a_{ij}}(-1)^s\begin{bmatrix}1-a_{ij}\\
s\end{bmatrix}_{q_i}E_i^{1-a_{ij}-s}E_jE_i^s=0,
$$
$$
\sum_{s=0}^{1-a_{ij}}(-1)^s\begin{bmatrix}1-a_{ij}\\
s\end{bmatrix}_{q_i}F_i^{1-a_{ij}-s}F_jF_i^s=0,
$$
with
$$q_i=q^{d_i},\qquad 1\le i\le l,$$
$$
\begin{bmatrix}m\\ n\end{bmatrix}_q=\frac{[m]_q!}{[n]_q![m-n]_q!},\qquad
[n]_q!=[n]_q\ldots[2]_q[1]_q,\qquad [n]_q=\frac{q^n-q^{-n}}{q-q^{-1}}.
$$

The Hopf algebra structure on $U_q$ is given by the following formulas:
$$
\triangle(E_i)=E_i\otimes 1+K_i\otimes E_i,\qquad\triangle(F_i)=F_i\otimes
K_i^{-1}+1\otimes F_i,\qquad\triangle(K_i)=K_i\otimes K_i,
$$
$$S(E_i)=-K_i^{-1}E_i,\qquad S(F_i)=-F_iK_i,\qquad S(K_i)=K_i^{-1},$$
$$\varepsilon(E_i)=\varepsilon(F_i)=0,\qquad \varepsilon(K_i)=1.$$

Any subset $\mathbb{L}\subset\{1,2,\ldots,l\}$ determines a Hopf subalgebra
$U_q\mathfrak{l}\subset U_q$ generated by
$$K_i^{\pm 1},\quad i=1,2,\ldots,l;\qquad E_j,F_j,\quad j\in\mathbb{L}.$$

We use the following notation:
$$
P=\mathbb{Z}^l,\qquad P_+=\mathbb{Z}_+^l,\qquad
P_+^\mathbb{L}=\{\boldsymbol{\lambda}=(\lambda_1,\lambda_2,\ldots,\lambda_l)
\in P|\:\lambda_i\ge 0\text{\ \ for\ \ }i\in\mathbb{L}\}.
$$
Let $V$ be a $U_q\mathfrak{l}$-module and
$$V_\mu=\{v\in V|\:K_iv=q_i^{\mu_i}v,\quad i=1,2,\ldots,l\}$$
the weight space, where $\mu=(\mu_1,\mu_2,\ldots,\mu_l)\in P$.

\begin{Definition}
We say $V$ is a {\bf weight} $U_q\mathfrak{l}$-module if
$V=\bigoplus\limits_{\mu\in P}V_\mu$.
\end{Definition}

\begin{Example}\label{Verma}
Consider the Verma module $U_q\mathfrak{l}$-module
$M(\mathfrak{l},\lambda)$ with the highest weight $\lambda\in
P_+^\mathbb{L}$. This is a module determined by its generator
$v(\mathfrak{l},\lambda)$ and the relations
$$
E_jv(\mathfrak{l},\lambda)=0,\quad j\in\mathbb{L};\qquad K_i^{\pm
1}v(\mathfrak{l},\lambda)=q_i^{\pm\lambda_i}v(\mathfrak{l},\lambda),\quad
i=1,2,\ldots,l.
$$
It is easy to prove that its submodule $K(\mathfrak{l},\lambda)$ generated
by $F_j^{\lambda_j+1}v(\mathfrak{l},\lambda)$, $j\in\mathbb{L}$, is the
unique submodule of finite codimension. In particular, the factor-module
$L(\mathfrak{l},\lambda)=M(\mathfrak{l},\lambda)/K(\mathfrak{l},\lambda)$
is simple (see \cite{Jant}).
\end{Example}

The following claims are well-known (see for instance \cite{Jant}). The
dimensions of the weight subspaces of the $U_q\mathfrak{l}$-modules
$L(\mathfrak{l},\lambda)$ are the same as in the classical case $q=1$. The
weight finite dimensional $U_q\mathfrak{l}$-modules are completely
reducible. Irreducible weight finite dimensional $U_q\mathfrak{l}$-modules
are isomorphic to $L(\mathfrak{l},\lambda)$ with $\lambda\in
P_+^\mathbb{L}$. The center of $U_q\mathfrak{l}$ (denoted by
$Z(U_q\mathfrak{l})$) admits an easy description via the Harish-Chandra
isomorphism \cite[p. 109]{Jant} and separates the simple weight finite
dimensional $U_q\mathfrak{l}$-modules \cite[p. 125]{Jant}. In other words,
for any $\lambda_1,\lambda_2\in P_+^\mathbb{L}$ there exists $z\in
Z(U_q\mathfrak{l})$ such that $z|_{L(\mathfrak{l},\lambda_1)}\ne
z|_{L(\mathfrak{l},\lambda_2)}$.

A module $V$ over $U_q\mathfrak{l}$ is called locally finite dimensional if
$\dim(U_q\mathfrak{l}\,v)<\infty$ for all $v\in V$. Introduce the notation
$C(\mathfrak{l},\mathfrak{l})_q$ for the full subcategory of weight locally
finite dimensional $U_q\mathfrak{l}$-modules.

It is easy to prove that any $U_q\mathfrak{l}$-module $V\in
C(\mathfrak{l},\mathfrak{l})_q$ is uniquely decomposable as a direct sum
$V=\bigoplus\limits_{\lambda\in P_+^\mathbb{L}}V_\lambda$ of its
submodules. Each $V_\lambda$ here is a multiple of
$L(\mathfrak{l},\lambda)$, $\lambda\in P_+^\mathbb{L}$. The sum is direct
since $Z(U_q\mathfrak{l})$ separates simple finite dimensional weight
$U_q\mathfrak{l}$-modules.

Let $\operatorname{End}V$ be the algebra of all endomorphisms of $V$
considered as a vector space. Now consider the $U_q\mathfrak{l}$-module
$L^\mathrm{univ}=\bigoplus\limits_{\lambda\in
P_+^\mathbb{L}}L(\mathfrak{l},\lambda)$ and the projections $P_\lambda$ in
$L^\mathrm{univ}$ onto the isotypic component $L(\mathfrak{l},\lambda)$
parallel to the sum of all other isotypic components. The algebra
$\operatorname{End}L^\mathrm{univ}$ is a $U_q\mathfrak{l}$-module
algebra.\footnote{If $\xi\in U_q\mathfrak{l}$ and
$\triangle\xi=\sum\limits_i\xi_i'\otimes\xi_i''$, then $(\xi
A)v=\sum\limits_i\xi_i'AS(\xi_i'')v$, $v\in L^\mathrm{univ}$,
$A\in\operatorname{End}L^\mathrm{univ}$, with $\triangle$ and $S$ being the
comultiplication and the antipode of $U_q\mathfrak{l}$.} The natural
homomorphism $U_q\mathfrak{l}\to\operatorname{End}L^\mathrm{univ}$ is
injective \cite[p. 76 -- 77]{Jant}; this allows one to identify
$U_q\mathfrak{l}$ with its image in $\operatorname{End}L^\mathrm{univ}$. We
are interested in the following $U_q\mathfrak{l}$-module subalgebras of
$\operatorname{End}L^\mathrm{univ}$:
$$
R(\mathfrak{l})_q\overset{\mathrm{def}}{=}\bigoplus\limits_{\lambda\in
P_+^\mathbb{L}}\operatorname{End}L(\mathfrak{l},\lambda),\qquad
F(\mathfrak{l})_q\overset{\mathrm{def}}{=}R(\mathfrak{l})_q\oplus
U_q\mathfrak{l}.\footnote{It is easy to prove that $R(\mathfrak{l})_q\cap
U_q\mathfrak{l}=0$. In fact, $U_q\mathfrak{l}$ has no zero divisors
\cite{DeConKac}, while the Harish-Chandra isomorphism allows one to produce
for any $a\in R(\mathfrak{l})_q$ a non-zero element $z\in
Z(U_q\mathfrak{l})$ such that $za=0$.}
$$
Obviously, $R(\mathfrak{l})_q$ is a two-sided ideal of $F(\mathfrak{l})_q$
generated by the projections $P_\lambda$.

$R(\mathfrak{l})_q$ is not a unital algebra, but it admits a distinguished
approximate identity. Specifically, any finite subset $\Lambda\subset
P_+^\mathbb{L}$ determines an idempotent
$\chi_\Lambda=\sum\limits_{\lambda\in\Lambda}P_\lambda$ in
$R(\mathfrak{l})_q$. Obviously,
$\chi_{\Lambda_1}\chi_{\Lambda_2}=\chi_{\Lambda_2}\chi_{\Lambda_1}=
\chi_{\Lambda_1\cap\Lambda_2}$, and for any $r\in R(\mathfrak{l})_q$ there
exists a finite subset $\Lambda\subset P_+^\mathbb{L}$ such that
$\chi_\Lambda r=r\chi_\Lambda=r$. With some abuse of terminology, one can
say that in $F(\mathfrak{l})_q$ we have
$$\lim_{\Lambda\uparrow P_+^\mathbb{L}}\chi_\Lambda=1.$$

A module $V$ over $R(\mathfrak{l})_q$ is called approximately unital if for
any $v\in V$ there exists a finite subset $\Lambda\subset P_+^\mathbb{L}$
such that $\chi_\Lambda v=v$.

We are going to demonstrate that the full subcategory of approximately
unital $R(\mathfrak{l})_q$-modules is canonically isomorphic to
$C(\mathfrak{l},\mathfrak{l})_q$.

\begin{Proposition}\label{iso_l}
(*). For any approximately unital $R(\mathfrak{l})_q$-module $V$ and any
$\xi\in F(\mathfrak{l})_q$, $v\in V$ there exists a limit
\begin{equation}\label{lim_l}
\xi v\overset{\mathrm{def}}{=}\lim_{\Lambda\uparrow
P_+^\mathbb{L}}(\xi\chi_\Lambda)v
\end{equation}
in the sense that $(\xi\chi_\Lambda)v$ does not depend on $\Lambda$ for
large $\Lambda$.

(**). The relation \eqref{lim_l} equips $V$ with a structure of
$F(\mathfrak{l})_q$-module and defines a functor from the category of
approximately unital $R(\mathfrak{l})_q$-modules into the category
$C(\mathfrak{l},\mathfrak{l})_q$.

(***). This functor is an isomorphism of categories.
\end{Proposition}

{\bf Proof.} Any approximately unital $R(\mathfrak{l})_q$-module $V$ admits
an embedding into an approximately unital $R(\mathfrak{l})_q$-module that
is a multiple of $L^\mathrm{univ}$. This observation implies the first two
statements, since they can be easily verified for $L^\mathrm{univ}$.

Any $U_q\mathfrak{l}$-module that is a multiple of $L^\mathrm{univ}$, is a
$F(\mathfrak{l})_q$-module. This allows one to construct an inverse functor
by embedding any $U_q\mathfrak{l}$-module $V\in
C(\mathfrak{l},\mathfrak{l})_q$ into a $U_q\mathfrak{l}$-module that is a
multiple of $L^\mathrm{univ}$. This implies (***). \hfill $\square$

\medskip

Note that the action of $P_\lambda\in R(\mathfrak{l})_q$ on vectors $v\in
V\in C(\mathfrak{l},\mathfrak{l})_q$ can be described without referring to
an embedding into a $U_q\mathfrak{l}$-module multiple to $L^\mathrm{univ}$.
More exactly, for any $\lambda\in P_+^\mathbb{L}$ there exists a sequence
$z_1(\lambda),z_2(\lambda),\ldots$ of elements of the center of
$U_q\mathfrak{l}$ such that $P_\lambda
v=\lim\limits_{n\to\infty}z_n(\lambda)v$ for all $v\in V\in
C(\mathfrak{l},\mathfrak{l})_q$. To see this, note that for any finite
subset $\lambda\in\Lambda\subset P_+^\mathbb{L}$ there exists such $z\in
Z(U_q\mathfrak{l})$ that $z|_{L(\mathfrak{l},\mu)}=0$ for
$\lambda\in\Lambda\setminus\{\lambda\}$ and
$z|_{L(\mathfrak{l},\lambda)}=1$. It remains to choose a sequence
$\Lambda_n\uparrow P_+^\mathbb{L}$.

\section{The category \boldmath $C(\mathfrak{g},\mathfrak{l})_q$ and the
algebra $R(\mathfrak{g},\mathfrak{l})_q$}

Let $\mathbb{G}\supset\mathbb{L}$ be a pair of subsets of
$\{1,2,\ldots,l\}$ and let $U_q\mathfrak{g}\supset U_q\mathfrak{l}$ be the
associated pair of Hopf subalgebras of $U_q$.

Consider the category of all $U_q\mathfrak{g}$-modules and its full
subcategory $C(\mathfrak{g},\mathfrak{l})_q$ that consists of weight
$U_q\mathfrak{l}$-locally finite dimensional $U_q\mathfrak{g}$-modules.

Let us turn to a description of the algebra
$F(\mathfrak{g},\mathfrak{l})_q$. This algebra is an important tool in
studying $U_q\mathfrak{g}$-modules from $C(\mathfrak{g},\mathfrak{l})_q$.
First, we consider a simple special case, namely $\mathbb{G}=\mathbb{L}$.
In this situation we define $F(\mathfrak{g},\mathfrak{l})_q$ as
$F(\mathfrak{l})_q$. We intend to construct
$F(\mathfrak{g},\mathfrak{l})_q$ along with embeddings
$U_q\mathfrak{g}\hookrightarrow F(\mathfrak{g},\mathfrak{l})_q$,
$F(\mathfrak{l})_q\hookrightarrow F(\mathfrak{g},\mathfrak{l})_q$, which
possess the following properties.
\begin{description}
\item[(F1)] The diagram
$$
\xymatrix{U_q\mathfrak{l}\ar@{^(->}[r]\ar@{_(->}[d] &
U_q\mathfrak{g}\ar@{_(->}[d]\\ F(\mathfrak{l})_q\ar@{^(->}[r] &
F(\mathfrak{g},\mathfrak{l})_q}
$$
commutes\footnote{In what follows $U_q\mathfrak{g}$ and $F(\mathfrak{l})_q$
are identified to their images under the embedding into
$F(\mathfrak{g},\mathfrak{l})_q$.}, and the subalgebras $U_q\mathfrak{g}$
and $F(\mathfrak{l})_q$ generate the algebra
$F(\mathfrak{g},\mathfrak{l})_q$.

\item[(F2)] The approximate identity $\{\chi_\Lambda\}$ of
$R(\mathfrak{l})_q$ is also an approximate identity of the two-sided ideal
$R(\mathfrak{g},\mathfrak{l})_q$ of $F(\mathfrak{g},\mathfrak{l})_q$,
generated by $R(\mathfrak{l})_q$.

\item[(F3)] For any approximately unital
$R(\mathfrak{g},\mathfrak{l})_q$-module $V$ and any $\xi\in
F(\mathfrak{g},\mathfrak{l})_q$, $v\in V$ there exists a limit
\begin{equation}\label{lim_g}
\xi v\overset{\mathrm{def}}{=}\lim_{\Lambda\uparrow
P_+^\mathbb{L}}(\xi\chi_\Lambda)v,
\end{equation}
and the equation \eqref{lim_g} equips $V$ with a structure of
$F(\mathfrak{g},\mathfrak{l})_q$-module.

The intersection of the kernels for all representations of
$F(\mathfrak{g},\mathfrak{l})_q$ obtained in this way is 0.

\item[(F4)] The functor from the category of approximately unital
$R(\mathfrak{g},\mathfrak{l})_q$-modules into
$C(\mathfrak{g},\mathfrak{l})_q$ given by \eqref{lim_g} is an isomorphism of
categories.
\end{description}

We are going to show that the algebra $F(\mathfrak{g},\mathfrak{l})_q$ and
the embeddings $U_q\mathfrak{g}\hookrightarrow
F(\mathfrak{g},\mathfrak{l})_q$, $F(\mathfrak{l})_q\hookrightarrow
F(\mathfrak{g},\mathfrak{l})_q$ are essentially uniquely determined by (F1)
-- (F4), i.e. there exists a unique isomorphism of algebras that respects
the embeddings.

Assign to any $U_q\mathfrak{l}$-module $V$ a $U_q\mathfrak{g}$-module
$P(V)=U_q\mathfrak{g}\otimes_{U_q\mathfrak{l}}V$. The following result is
well-known in the classical case $q=1$ (see \cite{KV}). It is proved in the
Appendix (see Corollary \ref{g_l}).

\begin{Lemma}\label{P_V}
1. If $V\in C(\mathfrak{l},\mathfrak{l})_q$ then $P(V)$ is a projective
object in $C(\mathfrak{g},\mathfrak{l})_q$.

2. For any $V\in C(\mathfrak{g},\mathfrak{l})_q$, the map
$$
P(V)\mapsto V,\qquad\xi\otimes v\mapsto\xi v,\qquad\xi\in
U_q\mathfrak{g},\quad v\in V,
$$
is a surjective morphism of $U_q\mathfrak{g}$-modules.
\end{Lemma}

The projective objects $P(V)$ of the category
$C(\mathfrak{g},\mathfrak{l})_q$ are called standard projective objects.

\medskip

Introduce the $U_q\mathfrak{g}$-module $V^\mathrm{univ}=P(L^\mathrm{univ})$.

\begin{Corollary}\label{projectors}
In the category of $U_q\mathfrak{l}$-modules one has
$$
V^\mathrm{univ}=\bigoplus_{\lambda\in
P_+^\mathbb{L}}V_\lambda^\mathrm{univ},
$$
where the $U_q\mathfrak{l}$-modules $V_\lambda^\mathrm{univ}$ are multiples
of $L(\mathfrak{l},\lambda)$.
\end{Corollary}

\begin{Corollary}\label{new}
For any module $M\in C(\mathfrak{g},\mathfrak{l})_q$ there exists a
surjective morphism $\widetilde{V}^\mathrm{univ}\to M$ in the category
$C(\mathfrak{g},\mathfrak{l})_q$, with $\widetilde{V}^\mathrm{univ}$ being
a multiple of the module $V^\mathrm{univ}$.
\end{Corollary}

Let us introduce the notation $\mathscr{P}_\lambda$ for the projection in
$V^\mathrm{univ}$ onto the $U_q\mathfrak{l}$-isotypic component
$V_\lambda^\mathrm{univ}$ parallel to the sum of all other
$U_q\mathfrak{l}$-isotypic components.

The category $C(\mathfrak{g},\mathfrak{l})_q$ and the category of
approximately unital $R(\mathfrak{g},\mathfrak{l})_q$-modules are closed
under the operations of direct sums, passage to submodules and
factor-modules. This allows one, using Lemma \ref{P_V} and Corollary
\ref{new}, to prove the following statement, which implies uniqueness of
$F(\mathfrak{g},\mathfrak{l})_q$ and the embeddings
$U_q\mathfrak{g}\hookrightarrow F(\mathfrak{g},\mathfrak{l})_q$,
$F(\mathfrak{l})_q\hookrightarrow F(\mathfrak{g},\mathfrak{l})_q$.

\begin{Proposition}\label{uniqueness}
Consider the algebra $F(\mathfrak{g},\mathfrak{l})_q$ and the embeddings
$U_q\mathfrak{g}\hookrightarrow F(\mathfrak{g},\mathfrak{l})_q$,
$F(\mathfrak{l})_q\hookrightarrow F(\mathfrak{g},\mathfrak{l})_q$ with the
properties (F1) -- (F4). Then the following statements hold.

1. The representation of $U_q\mathfrak{g}$ in $V^\mathrm{univ}$ determines
an approximately unital representation of $R(\mathfrak{g},\mathfrak{l})_q$
and a faithful representation of $F(\mathfrak{g},\mathfrak{l})_q$ in
$V^\mathrm{univ}$.

2. The image of $F(\mathfrak{g},\mathfrak{l})_q$ under the embedding into
$\operatorname{End}V^\mathrm{univ}$ is the subalgebra of
$\operatorname{End}V^\mathrm{univ}$ generated by elements of
$U_q\mathfrak{g}\subset\operatorname{End}V^\mathrm{univ}$ and projections
$\mathscr{P}_\lambda$, $\lambda\in P_+^\mathbb{L}$.
\end{Proposition}

\medskip

Now we return to the construction of the algebra
$F(\mathfrak{g},\mathfrak{l})_q$. In what follows we will identify the
algebras $U_q\mathfrak{g}$, $F(\mathfrak{l})_q$ with their images under the
embeddings into $\operatorname{End}V^\mathrm{univ}$.\footnote{The
representation of $U_q\mathfrak{g}$ in $V^\mathrm{univ}$ is faithful since
for any weight finite dimensional $U_q\mathfrak{g}$-module $W$ there exists
a surjective morphism of $U_q\mathfrak{g}$-modules $V^\mathrm{univ}\to W$
(Corollary \ref{new}). The representation of $F(\mathfrak{l})_q$ in
$V^\mathrm{univ}$ is faithful since $L^\mathrm{univ}\hookrightarrow
V^\mathrm{univ}$.} In particular,
$$
P_\lambda\mapsto\mathscr{P}_\lambda,\quad \lambda\in P_+^\mathbb{L}.
$$

The following auxiliary statement is proved in the Appendix.

\begin{Lemma}\label{hard_lemma}
Given any $\xi\in U_q\mathfrak{g}$, $\lambda\in P_+^\mathbb{L}$, there
exists a finite subset $\Lambda\subset P_+^\mathbb{L}$ such that
$$
\xi\mathscr{P}_\lambda=\chi_\Lambda\xi\mathscr{P}_\lambda,\qquad
\mathscr{P}_\lambda\xi=\mathscr{P}_\lambda\xi\chi_\Lambda.
$$
\end{Lemma}

\medskip

Consider the smallest subalgebra
$F(\mathfrak{g},\mathfrak{l})_q\subset\operatorname{End}V^\mathrm{univ}$
that contains all the elements of $U_q\mathfrak{g}$ and all the projections
$\mathscr{P}_\lambda$, $\lambda\in P_+^\mathbb{L}$.

\begin{Theorem}\label{existence}
The algebra $F(\mathfrak{g},\mathfrak{l})_q$ along with the embeddings
$U_q\mathfrak{g}\hookrightarrow F(\mathfrak{g},\mathfrak{l})_q$,
$F(\mathfrak{l})_q\hookrightarrow F(\mathfrak{g},\mathfrak{l})_q$, satisfy
(F1) -- (F4).
\end{Theorem}

{\bf Proof.} (F1) is obviously satisfied.

Let us prove that $\{\chi_\Lambda\}$ is a left approximate identity in
$R(\mathfrak{g},\mathfrak{l})_q$. Every $a\in
R(\mathfrak{g},\mathfrak{l})_q$ has the form
$a=\sum\limits_{i=1}^{N(a)}\xi_i\mathscr{P}_{\lambda_i}\eta_i$, with
$\xi_i\in U_q\mathfrak{g}$, $\eta_i\in F(\mathfrak{g},\mathfrak{l})_q$,
$\lambda_i\in P_+^\mathbb{L}$. It follows from Lemma \ref{hard_lemma} that
$\chi_\Lambda\xi_i\mathscr{P}_{\lambda_i}=\xi_i\mathscr{P}_{\lambda_i}$ for
some finite subset $\Lambda\subset P_+^\mathbb{L}$ and all
$i=1,2,\ldots,N(a)$. Hence $\chi_\Lambda a=a$. One can prove in a similar
way that $\{\chi_\Lambda\}$ is a right approximate identity in
$R(\mathfrak{g},\mathfrak{l})_q$. Hence (F2) is satisfied.

Consider an approximately unital $R(\mathfrak{g},\mathfrak{l})_q$-module
$V$. Given any vector $v\in V$, there exists such finite subset
$\Lambda\subset P_+^\mathbb{L}$ that $\chi_\Lambda v=v$. Thus for all
$\Lambda'\supset\Lambda$ we have
$$
(\xi\chi_{\Lambda'})v=(\xi\chi_{\Lambda'})\chi_\Lambda
v=(\xi\chi_{\Lambda'}\chi_\Lambda)v=(\xi\chi_\Lambda)v,
$$
and hence the element
$$
\xi v\overset{\mathrm{def}}{=}\lim_{\Lambda\uparrow
P_+^\mathbb{L}}(\xi\chi_\Lambda)v
$$
is well defined. It is easy to prove that $(\xi\eta)v=\xi(\eta v)$ for any
$\xi,\eta\in F(\mathfrak{g},\mathfrak{l})_q$, $v\in V$. In fact, there
exist such finite subsets $\Lambda',\Lambda''\subset P_+^\mathbb{L}$ that
$$
\chi_{\Lambda'}v=v,\qquad
\chi_{\Lambda''}\eta\chi_{\Lambda'}=\eta\chi_{\Lambda'},
$$
because $\eta\chi_{\Lambda'}\in R(\mathfrak{g},\mathfrak{l})_q$ and
$\{\chi_\Lambda\}$ is an approximate identity. Hence,
\begin{multline*}
(\xi\eta)v=\lim_{\Lambda\uparrow
P_+^\mathbb{L}}(\xi\eta\chi_\Lambda)v=\xi(\eta\chi_{\Lambda'})v=
(\xi\chi_{\Lambda''})(\eta\chi_{\Lambda'})v=\lim_{\Lambda\uparrow
P_+^\mathbb{L}}(\xi\chi_\Lambda)(\eta\chi_{\Lambda'})v=
\\ =\xi(\eta\chi_{\Lambda'}v)=\xi(\eta v).
\end{multline*}
It follows that $V$ is an $F(\mathfrak{g},\mathfrak{l})_q$-module, and the
the first of the conditions (F3) is satisfied. The second condition, which
requires the existence of a faithful representation as above, obviously
holds.

The arguments below show that the $F(\mathfrak{g},\mathfrak{l})_q$-modules
in question are locally $U_q\mathfrak{l}$-finite dimensional and weight:
$$
\dim(U_q\mathfrak{l}\,v)=\dim
U_q\mathfrak{l}(\chi_{\Lambda'}v)=\dim(U_q\mathfrak{l}\,\chi_{\Lambda'})v\le
\dim R(\mathfrak{l})_qv<\infty.
$$
Thus we have constructed a functor from the category of approximately
unital $R(\mathfrak{g},\mathfrak{l})_q$-modules into
$C(\mathfrak{g},\mathfrak{l})_q$. What remains is to construct an inverse
functor.

Suppose $V\in C(\mathfrak{g},\mathfrak{l})_q$ and let
$V=\bigoplus\limits_{\mu\in P_+^\mathbb{L}}V_\mu$ be a decomposition of $V$
into a sum of $U_q\mathfrak{l}$-isotypic components. Let $\pi$ be the
representation of $U_q\mathfrak{g}$ in $V$, and denote by $\Pi_\lambda$ the
projection in $V$ onto the isotypic component of $V_\lambda$ parallel
to $\bigoplus\limits_{\substack{\mu\in P_+^\mathbb{L}\\
\mu\ne\lambda}}V_\mu$.

It suffices to prove existence and uniqueness of an extension
$\widetilde{\pi}$ of $\pi$ to $F(\mathfrak{g},\mathfrak{l})_q$ that
possesses the following properties:

i) $\widetilde{\pi}(\mathscr{P}_\lambda)=\Pi_\lambda$, $\lambda\in
P_+^\mathbb{L}$,

ii) $\widetilde{\pi}|_{R(\mathfrak{g},\mathfrak{l})_q}$ is an approximately
unital representation of $R(\mathfrak{g},\mathfrak{l})_q$,

iii) $\pi(\xi)v=\lim\limits_{\Lambda\uparrow
P_+^\mathbb{L}}\widetilde{\pi}(\xi\chi_\Lambda)v$, $v\in V$, $\xi\in
U_q\mathfrak{g}$.

Uniqueness of such extension is obvious. To prove its existence, we consider
subsequently the following cases:
\begin{enumerate}
\item $V=V^\mathrm{univ}$;

\item $V=\widetilde{V}^\mathrm{univ}$, with $\widetilde{V}^\mathrm{univ}$
being a multiple of the $U_q\mathfrak{g}$-module $V^\mathrm{univ}$;

\item $V$ is a submodule of $\widetilde{V}^\mathrm{univ}$;

\item $V$ is a standard projective object in
$C(\mathfrak{g},\mathfrak{l})_q$;

\item (the general case) $V\in C(\mathfrak{g},\mathfrak{l})_q$.
\end{enumerate}

In the case 1) the desired statement follows from the definition of
$R(\mathfrak{g},\mathfrak{l})_q$. A passage from 1) to 2) and to 3) is
obvious. To pass from 3) to 4), it is sufficient to use the fact that every
$U_q\mathfrak{g}$-module $P(V)$, $V\in C(\mathfrak{l},\mathfrak{l})_q$,
admits an embedding into a multiple $\widetilde{V}^\mathrm{univ}$ of the
$U_q\mathfrak{g}$-module $V^\mathrm{univ}$. This is because every locally
finite dimensional $U_q\mathfrak{l}$-module admits an embedding into a
$U_q\mathfrak{l}$-module $\widetilde{L}^\mathrm{univ}$, which is a multiple
of $L^\mathrm{univ}$ and has the property
$P(\widetilde{L}^\mathrm{univ})=\widetilde{V}^\mathrm{univ}$. Finally, to
pass from 4) to 5), one can apply the existence of a surjective morphism
$P(V)\to V$ in the category $C(\mathfrak{g},\mathfrak{l})_q$. \hfill
$\square$

\medskip

Recall that the Weyl group $W$ is generated by simple reflections
$s_1,s_2,\ldots,s_l$ and acts on the weight lattice $P$.

\begin{Proposition}\label{ext}
1. There exists a unique one-dimensional representation
$\widetilde{\varepsilon}$ of $F(\mathfrak{g},\mathfrak{l})_q$ with
$$
\widetilde{\varepsilon}(\xi)=\varepsilon(\xi),\quad\xi\in
U_q\mathfrak{g};\qquad\widetilde{\varepsilon}(P_\lambda)=
\begin{cases}
1, & \lambda=0,
\\ 0, & \lambda\in P_+^\mathbb{L}\setminus\{0\}.
\end{cases}
$$

2. There exists a unique anti-automorphism $\widetilde{S}$ of the algebra
$F(\mathfrak{g},\mathfrak{l})_q$ such that
$$
\widetilde{S}(\xi)=S(\xi),\quad\xi\in
U_q\mathfrak{g};\qquad\widetilde{S}(P_\lambda)=P_{-w_0^\mathbb{L}\lambda},
\quad\lambda\in P_+^\mathbb{L},
$$
where $w_0^\mathbb{L}$ is the longest element of the Weyl group
$W_\mathbb{L}\subset W$ generated by the simple reflections $s_i$,
$i\in\mathbb{L}$.
\end{Proposition}

{\bf Proof.} Uniqueness of $\widetilde{\varepsilon}$ and $\widetilde{S}$ is
obvious. Existence of $\widetilde{\varepsilon}$ has been demonstrated in
the proof of Theorem \ref{existence}. Turn to proving existence of
$\widetilde{S}$. In the special case $\mathbb{G}=\mathbb{L}$ it follows
from the definition of the algebra
$F(\mathfrak{l})_q\hookrightarrow\operatorname{End}L^\mathrm{univ}$ and the
well-known isomorphism \cite[p. 168]{Bou7-8}
$$
L(\mathfrak{l},\lambda)^*\overset{\sim}{\to}
L(\mathfrak{l},-w_0^\mathbb{L}\lambda),\qquad v(\mathfrak{l},\lambda)\mapsto
v(\mathfrak{l},w_0^\mathbb{L}\lambda).
$$

To pass from the special case $\mathbb{G}=\mathbb{L}$ to the general case
$\mathbb{G}\supset\mathbb{L}$, consider the algebra
$F(\mathfrak{g},\mathfrak{l})_q^\mathrm{op}$ which is derived from
$F(\mathfrak{g},\mathfrak{l})_q$ by replacement of the multiplication by
the opposite one. $U_q\mathfrak{g}$ can be embedded into
$F(\mathfrak{g},\mathfrak{l})_q^\mathrm{op}$:
$$
U_q\mathfrak{g}\hookrightarrow
F(\mathfrak{g},\mathfrak{l})_q^\mathrm{op},\qquad\xi\mapsto S(\xi).
$$
The algebra $F(\mathfrak{l})_q$ also admits an embedding into
$F(\mathfrak{g},\mathfrak{l})_q^\mathrm{op}$:
$$
U_q\mathfrak{l}\hookrightarrow
F(\mathfrak{g},\mathfrak{l})_q^\mathrm{op},\quad\xi\mapsto S(\xi);\qquad
P_\lambda\mapsto P_{-w_0^\mathbb{L}\lambda}.
$$
It is easy to verify (F1) -- (F4) for this pair of embeddings. What remains
is to use the uniqueness of such pair, Proposition \ref{uniqueness}. \hfill
$\square$

\medskip

The category $C(\mathfrak{g},\mathfrak{l})_q$ is canonically isomorphic to
the category of approximately unital
$R(\mathfrak{g},\mathfrak{l})_q$-modules. In what follows we identify these
categories.

Consider the Hopf subalgebra $U_q\mathfrak{q}_\mathbb{L}^+\subset
U_q\mathfrak{g}$ generated by
$$
E_i,\quad i\in\mathbb{G};\qquad F_j,\quad j\in\mathbb{L};\qquad K_m^{\pm
1},\quad m=1,2,\ldots,l.
$$
and the Hopf subalgebra $U_q\mathfrak{q}_\mathbb{L}^-\subset
U_q\mathfrak{g}$ generated by
$$
E_i,\quad i\in\mathbb{L};\qquad F_j,\quad j\in\mathbb{G};\qquad K_m^{\pm
1},\quad m=1,2,\ldots,l.
$$

The following statement is well known in the classical case $q=1$ (see
\cite[p. 90]{KV}).

\begin{Proposition}\label{like_KV}
The linear maps
\begin{equation}\label{first}
U_q\mathfrak{g}\otimes_{U_q\mathfrak{l}}R(\mathfrak{l})_q\to
R(\mathfrak{g},\mathfrak{l})_q,\qquad\xi\otimes r\mapsto\xi r,
\end{equation}
$$
R(\mathfrak{l})_q\otimes_{U_q\mathfrak{l}}U_q\mathfrak{g}\to
R(\mathfrak{g},\mathfrak{l})_q,\qquad r\otimes\xi\mapsto r\xi,
$$
are bijective.
\end{Proposition}

{\bf Proof.} By Proposition \ref{ext}, it suffices to consider the linear
map \eqref{first}. Let us prove that it is onto. Note that for any
$\lambda\in P_+^\mathbb{L}$ and any finite subset $\Lambda\subset
P_+^\mathbb{L}$ containing $\lambda$ there exists an element $z$ in the
center $Z(U_q\mathfrak{l})$ of the algebra $U_q\mathfrak{l}$ such that
$P_\lambda=z\chi_\Lambda$. On the other hand, for any $\lambda\in
P_+^\mathbb{L}$, $\xi\in U_q\mathfrak{g}$ there exist finite subsets
$\Lambda',\Lambda''\ni\lambda$ such that
$$
P_\lambda\xi=P_\lambda\xi\chi_{\Lambda'},\qquad
\xi\chi_{\Lambda'}=\chi_{\Lambda''}\xi\chi_{\Lambda'}.
$$
Hence $P_\lambda\xi=(P_\lambda\chi_{\Lambda''})\xi\chi_{\Lambda'}=
z\chi_{\Lambda''}\xi\chi_{\Lambda'}=z\xi\chi_{\Lambda'}$. Thus for any
$\xi\in U_q\mathfrak{g}$ and any finite subset $\Lambda\subset
P_+^\mathbb{L}$ there exist such $\widetilde{\xi}\in U_q\mathfrak{g}$ and a
finite subset $\widetilde{\Lambda}\subset P_+^\mathbb{L}$ that
$\chi_\Lambda\xi=\widetilde{\xi}\chi_{\widetilde{\Lambda}}$. This can be
readily used to prove that the linear map \eqref{first} is onto.

Prove that it is injective. One can use well-known results on bases in
quantum universal enveloping algebras \cite[Chapter 8]{Jant} to derive a
decomposition
$U_q\mathfrak{g}=U_q\mathfrak{q}_\mathbb{L}^-\otimes_{U_q\mathfrak{l}}
U_q\mathfrak{q}_\mathbb{L}^+$. It follows from this decomposition and Lemma
\ref{finite_rk} that $U_q\mathfrak{g}$ is a free right
$U_q\mathfrak{l}$-module. Choose a free basis $\{e_i\}$. Let $a$ be a
non-zero element of
$U_q\mathfrak{g}\otimes_{U_q\mathfrak{l}}R(\mathfrak{l})_q$. It has the
form
$$a=\sum_ie_i\otimes r_i,\qquad r_i\in R(\mathfrak{l})_q.$$
Let $\widetilde{a}$ be the image of $a$ under the map \eqref{first}. We
treat it as an element of $\operatorname{End}V^\mathrm{univ}$. It suffices
to prove that the restriction of the linear map $\widetilde{a}$ to the
subspace $\{1\otimes v|\:v\in L^\mathrm{univ}\}$ is non-zero. It is easy to
verify that $\widetilde{a}(1\otimes v)=\sum\limits_ie_i\otimes r_iv$ for
any $v\in L^\mathrm{univ}$. It remains to use our choice of $\{e_i\}$ and
the fact that the representation of $U_q\mathfrak{l}$ in $L^\mathrm{univ}$
is faithful. \hfill $\square$

\medskip

\section{The functors \boldmath $\mathrm{ind}$ and $\Pi$}

Consider two pairs of subsets
$$
\mathbb{L}\subset\mathbb{G}\subset\{1,2,\ldots,l\},\qquad
\mathbb{L}_1\subset\mathbb{G}_1\subset\{1,2,\ldots,l\},
$$
with $\mathbb{L}_1\subset\mathbb{L}$, $\mathbb{G}_1\subset\mathbb{G}$.
Obviously, one has embeddings of the associated Hopf subalgebras
$$
U_q\mathfrak{l}_1\subset U_q\mathfrak{l},\qquad U_q\mathfrak{g}_1\subset
U_q\mathfrak{g}.
$$

$R(\mathfrak{g},\mathfrak{l})_q$ is a left ideal in
$F(\mathfrak{g},\mathfrak{l})_q$ and hence is a left
$U_q\mathfrak{g}_1$-module. We claim that $R(\mathfrak{g},\mathfrak{l})_q$
is an approximately unital left
$R(\mathfrak{g}_1,\mathfrak{l}_1)_q$-module. In fact,
$R(\mathfrak{g},\mathfrak{l})_q$ is an approximately unital left
$R(\mathfrak{g},\mathfrak{l})_q$-module, hence a $U_q\mathfrak{g}$-module
of the category $C(\mathfrak{g},\mathfrak{l})_q$ given by \eqref{lim_l}.
Since $U_q\mathfrak{l}_1\subset U_q\mathfrak{l}$, the module
$R(\mathfrak{g},\mathfrak{l})_q$ is also in the category
$C(\mathfrak{g},\mathfrak{l}_1)_q$. Finally, in view of
$U_q\mathfrak{g}_1\subset U_q\mathfrak{g}$, we are inside the category
$C(\mathfrak{g}_1,\mathfrak{l}_1)_q$, which is equivalent to our claim. In
a similar way, one can prove that $R(\mathfrak{g},\mathfrak{l})_q$ is a
right $R(\mathfrak{g}_1,\mathfrak{l}_1)_q$-module.\footnote{One can see
from the proof that $R(\mathfrak{g},\mathfrak{l})_q$ is a
$R(\mathfrak{g}_1,\mathfrak{l}_1)_q$-bimodule.}

Throughout the rest of this section we follow \cite{KV} and replace the
groups involved therein with the corresponding quantum universal enveloping
algebras. Introduce the functor
$P_{\mathfrak{g}_1,\mathfrak{l}_1}^{\mathfrak{g},\mathfrak{l}}$ from the
category $C(\mathfrak{g}_1,\mathfrak{l}_1)_q$ to the category
$C(\mathfrak{g},\mathfrak{l})_q$ by defining it on objects from
$C(\mathfrak{g}_1,\mathfrak{l}_1)_q$ as follows:
$$
P_{\mathfrak{g}_1,\mathfrak{l}_1}^{\mathfrak{g},\mathfrak{l}}(Z)=
R(\mathfrak{g},\mathfrak{l})_q\otimes_{R(\mathfrak{g}_1,\mathfrak{l}_1)_q}Z.
$$
The action on morphisms is defined in an obvious way.

$C(\mathfrak{g}_1,\mathfrak{l}_1)_q$ has enough projectives and the functor
$P_{\mathfrak{g}_1,\mathfrak{l}_1}^{\mathfrak{g},\mathfrak{l}}$ is
covariant and right exact.\footnote{See \cite[p. 840]{KV}} Hence one has
well defined derived functors: $\left(P_{\mathfrak{g}_1,\mathfrak{l}_1}
^{\mathfrak{g},\mathfrak{l}}\right)_j$, $j\in\mathbb{Z}_+$, from
$C(\mathfrak{g}_1,\mathfrak{l}_1)_q$ to $C(\mathfrak{g},\mathfrak{l})_q$.

Consider the two special cases: $\mathfrak{l}_1=\mathfrak{l}$ and
$\mathfrak{g}_1=\mathfrak{g}$. Start from the first one. Let
$\operatorname{ind}_{\mathfrak{g}_1,\mathfrak{l}}
^{\mathfrak{g},\mathfrak{l}}$ be the functor from the category
$C(\mathfrak{g}_1,\mathfrak{l})_q$ to the category
$C(\mathfrak{g},\mathfrak{l})_q$ defined on objects as follows:
$$
\operatorname{ind}_{\mathfrak{g}_1,\mathfrak{l}}^{\mathfrak{g},\mathfrak{l}}
(Z)=U_q\mathfrak{g}\otimes_{U_q\mathfrak{g}_1}Z.
$$
The action on morphisms is defined in an obvious way. Just as in the
classical case $q=1$, one gets an isomorphism of functors
$P_{\mathfrak{g}_1,\mathfrak{l}}^{\mathfrak{g},\mathfrak{l}}$ and
$\operatorname{ind}_{\mathfrak{g}_1,\mathfrak{l}}
^{\mathfrak{g},\mathfrak{l}}$. Describe briefly this construction. The
functor $\operatorname{ind}_{\mathfrak{g}_1,\mathfrak{l}}
^{\mathfrak{g},\mathfrak{l}}$ is left adjoint to the forgetful functor
$\mathscr{F}_{\mathfrak{g},\mathfrak{l}}
^{\mathfrak{g}_1,\mathfrak{l}}:C(\mathfrak{g},\mathfrak{l})_q\to
C(\mathfrak{g}_1,\mathfrak{l})_q$ determined by the embedding
$U_q\mathfrak{g}_1\to U_q\mathfrak{g}$. On the other hand,
$P_{\mathfrak{g}_1,\mathfrak{l}}^{\mathfrak{g},\mathfrak{l}}$ is left
adjoint to the functor $(\mathscr{F}^\vee)
_{\mathfrak{g},\mathfrak{l}}^{\mathfrak{g}_1,\mathfrak{l}}$, to be defined
below (cf. \cite[Proposition 2.34]{KV}). The functor $(\mathscr{F}^\vee)
_{\mathfrak{g},\mathfrak{l}}^{\mathfrak{g}_1,\mathfrak{l}}:
C(\mathfrak{g},\mathfrak{l})_q\to C(\mathfrak{g}_1,\mathfrak{l})_q$ is
defined on objects of the category $C(\mathfrak{g},\mathfrak{l})_q$ as
follows:
$$
(\mathscr{F}^\vee)_{\mathfrak{g},\mathfrak{l}}
^{\mathfrak{g}_1,\mathfrak{l}}(X)=
\operatorname{Hom}_{R(\mathfrak{g},\mathfrak{l})_q}
(R(\mathfrak{g},\mathfrak{l})_q,X)_\mathfrak{l}.
$$
A structure of $U_q\mathfrak{g}_1$-module in
$\operatorname{Hom}_{R(\mathfrak{g},\mathfrak{l})_q}
(R(\mathfrak{g},\mathfrak{l})_q,X)$ is imposed via the structure of right
$U_q\mathfrak{g}_1$-module in $R(\mathfrak{g},\mathfrak{l})_q$, and the
subscript $\mathfrak{l}$ stands for distinguishing the maximal submodule in
the category $C(\mathfrak{g}_1,\mathfrak{l})_q$. The action of
$(\mathscr{F}^\vee)_{\mathfrak{g},\mathfrak{l}}
^{\mathfrak{g}_1,\mathfrak{l}}$ on morphisms is defined in an obvious way.
What remains is to construct an isomorphism of functors
$\mathscr{F}_{\mathfrak{g},\mathfrak{l}}^{\mathfrak{g}_1,\mathfrak{l}}
\overset{\sim}{\to}(\mathscr{F}^\vee)_{\mathfrak{g},\mathfrak{l}}^
{\mathfrak{g}_1,\mathfrak{l}}$ (cf. \cite[Proposition 2.33]{KV}). Let $X\in
C(\mathfrak{g},\mathfrak{l})_q$ and $x\in X$. Associate to every $x$ a
morphism of $R(\mathfrak{g},\mathfrak{l})_q$-modules given by
$$R(\mathfrak{g},\mathfrak{l})_q\to X,\qquad r\mapsto rx.$$
One can verify that the map
$X\to\operatorname{Hom}_{R(\mathfrak{g},\mathfrak{l})_q}
(R(\mathfrak{g},\mathfrak{l})_q,X)$, which arises this way, is a morphism
of $U_q\mathfrak{g}_1$-modules and provides the desired isomorphism of
functors.

An important consequence is the observation that for any $V\in
C(\mathfrak{l},\mathfrak{l})_q$ the standard projective object $P(V)$ in
the category $C(\mathfrak{g},\mathfrak{l})_q$ is canonically isomorphic to
$P_{\mathfrak{l},\mathfrak{l}}^{\mathfrak{g},\mathfrak{l}}(V)$.

In the second special case $\mathfrak{g}_1=\mathfrak{g}$ the functor in
question is called the Bernstein functor and it is denoted by
$\Pi_{\mathfrak{g},\mathfrak{l}_1}^{\mathfrak{g},\mathfrak{l}}$:
$$
\Pi\equiv\Pi_{\mathfrak{g},\mathfrak{l}_1}^{\mathfrak{g},\mathfrak{l}}(Z)=
R(\mathfrak{g},\mathfrak{l})_q\otimes_{R(\mathfrak{g},\mathfrak{l}_1)_q}Z,
\qquad Z\in C(\mathfrak{g}_1,\mathfrak{l}_1)_q.
$$

In the classical case $q=1$ the derived functors $\Pi_j$ are crucial in
constructing unitarizable Harish-Chandra modules via the Vogan-Zuckerman
cohomological induction \cite{KV}. Turn to describing a quantum analogue
for this method.

\section{A quantum analogue for cohomological induction}

We assume in the sequel that $\mathbb{G}=\{1,2,\ldots,l\}$, and hence
$U_q\mathfrak{g}=U_q$.

Recall (see \cite{Helg, Kor}) a definition of the pairs
$(\mathfrak{g},\mathfrak{k})$ used in a construction of bounded symmetric
domains via the Harish-Chandra embedding. Here $\mathfrak{g}$ is a simple
complex Lie algebra of type $A$, $B$, $C$, $D$, $E_6$, $E_7$. A Hopf
subalgebra $U_q\mathfrak{k}\subset U_q\mathfrak{g}$ is generated by $E_i$,
$F_i$, $i\in\mathbb{K}$; $K_j^{\pm 1}$, $j=1,2,\ldots,l$, with
$\mathbb{K}=\{1,2,\ldots,l\}\setminus\{l_0\}$. We assume also that the
simple root $\alpha_{l_0}$ has coefficient $1$ in the decomposition of the
maximal root
\begin{equation}\label{herm_l0}
\delta=\sum_{i=1}^ln_i\alpha_i,\qquad n_{l_0}=1.
\end{equation}
In what follows the notation $(\mathfrak{g},\mathfrak{k})$ will stand only
for such pairs.

To clarify the precise nature of condition \eqref{herm_l0}, equip the Lie
algebra $\mathfrak{g}$ with a grading as follows: $\deg E_{l_0}=1$, $\deg
F_{l_0}=-1$, $\deg H_{l_0}=0$, $\deg E_j=\deg F_j=\deg H_j=0$ for $j\ne
l_0$. The automorphism $\theta$ of $\mathfrak{g}$ given by the formula
$$\theta(\xi)=(-1)^{\deg\xi}\xi$$
is an involution. Now \eqref{herm_l0} implies that
$$
\mathfrak{g}=\mathfrak{p}^-\oplus\mathfrak{k}\oplus\mathfrak{p}^+,\qquad
\mathfrak{p}^\pm=\{\xi\in\mathfrak{g}|\:\deg\xi=\pm 1\}.
$$

It is worthwhile to note that the pairs $(\mathfrak{g},\mathfrak{k})$ in
question are complexifications of the pairs
$(\mathfrak{g}_0,\mathfrak{k}_0)$, with $\mathfrak{g}_0$ being the Lie
subalgebra of the automorphism group of an irreducible bounded symmetric
domain, and $\mathfrak{k}_0$ being the Lie algebra of the stabilizer of a
point in this domain \cite{Helg}.

Equip $U_q\mathfrak{g}$ with a structure of Hopf $*$-algebra via the
involution $*$ defined as follows:
$$
E_j^*=
\begin{cases}
-K_jF_j, & j=l_0
\\ K_jF_j, & j\ne l_0
\end{cases};\qquad
F_j^*=
\begin{cases}
-E_jK_j^{-1}, & j=l_0
\\ E_jK_j^{-1}, & j\ne l_0
\end{cases};
$$
$$(K_j^{\pm 1})^*=K_j^{\pm 1},\qquad j=1,2,\ldots,l.$$
A module $V$ over a Hopf $*$-algebra $A$ is said to be unitarizable if it
admits a positive definite invariant form $(\cdot,\cdot)$:
$$(av_1,v_2)=(v_1,a^*v_2),\qquad a\in A,\;v_1,v_2\in V.$$

The cohomological induction is among the tools for constructing
unitarizable modules of the category $C(\mathfrak{g},\mathfrak{k})$ in the
case $q=1$. Describe a q-analogue for this method.

Let $\mathbb{L}\subset\{1,2,\ldots,l\}$, $\mathbb{L}\not\subset\mathbb{K}$,
and $U_q\mathfrak{l}$ is the Hopf subalgebra corresponding to the subset
$\mathbb{L}$. Obviously, $U_q\mathfrak{l}$ inherits the structure of Hopf
$*$-algebra. We intend to derive a unitarizable $U_q\mathfrak{g}$-module of
the category $C(\mathfrak{g},\mathfrak{k})_q$ starting from a unitarizable
$U_q\mathfrak{l}$-module of the category
$C(\mathfrak{l},\mathfrak{l}\cap\mathfrak{k})_q$.
\footnote{$U_q(\mathfrak{l}\cap\mathfrak{k})$ is a Hopf subalgebra
corresponding to the subset $\mathbb{L}\cap\mathbb{K}$.}

If $\lambda=(\lambda_1,\lambda_2,\ldots,\lambda_l)\in P$ and $\lambda_j=0$
for all $j\in\mathbb{L}$, one has the following well defined
one-dimensional $U_q\mathfrak{l}$-module $\mathbb{C}_\lambda$:
$$
E_j\mathbf{1}=F_j\mathbf{1}=0,\quad j\in\mathbb{L};\qquad K_j^{\pm
1}\mathbf{1}=q_j^{\pm\lambda_j}\mathbf{1},\quad j=1,2,\ldots,l.
$$
As an important example of such linear functional one should mention the
difference $\rho_\mathfrak{u}=\rho-\rho_\mathfrak{l}$ between the half-sum
$\rho$ of positive roots of the Lie algebra $\mathfrak{g}$ and the half-sum
$\rho_\mathfrak{l}$ of positive roots of the Lie algebra $\mathfrak{l}$:
$$
\rho_{\mathfrak{u}}(H_j)=
\begin{cases}
1, & j\not\in\mathbb{L},
\\ 0, & j\in\mathbb{L}.
\end{cases}
$$

Every $U_q\mathfrak{l}$-module $Z\in
C(\mathfrak{l},\mathfrak{l}\cap\mathfrak{k})_q$ determines a
$U_q\mathfrak{l}$-module $Z^\#=Z\otimes\mathbb{C}_{2\rho_\mathfrak{u}}\in
C(\mathfrak{l},\mathfrak{l}\cap\mathfrak{k})_q$. Equip it with a structure
of $U_q\mathfrak{q}_\mathbb{L}^-$-module via the surjective morphism of
algebras $U_q\mathfrak{q}_\mathbb{L}^-\to U_q\mathfrak{l}$:
$$
F_i\mapsto
\begin{cases}
F_i, & i\in\mathbb{L}
\\ 0, & i\in\{1,2,\ldots,l\}\setminus\mathbb{L}
\end{cases},\quad
E_j\mapsto E_j,\quad K_j^{\pm 1}\mapsto K_j^{\pm 1},\;j\in\{1,2,\ldots,l\}.
$$
The generalized Verma module
$\mathrm{ind}_{\mathfrak{q}_\mathbb{L}^-}^\mathfrak{g}Z^\#
\overset{\mathrm{def}}{=}U_q\mathfrak{g}
\otimes_{U_q\mathfrak{q}_\mathbb{L}^-}Z^\# $ belongs to the category
$C(\mathfrak{g},\mathfrak{l})_q$ (see Lemma \ref{finite_rk}). Hence the
equality
$$
\mathscr{L}_j(Z)=\left(\Pi_{\mathfrak{g},\mathfrak{l}\cap\mathfrak{k}}
^{\mathfrak{g},\mathfrak{k}}\right)_j
\left(\mathrm{ind}_{\mathfrak{q}_\mathbb{L}^-}^\mathfrak{g}(Z^\#)\right)
$$
determines functors $\mathscr{L}_j$, $j\in\mathbb{Z}_+$, from
$C(\mathfrak{l},\mathfrak{l}\cap\mathfrak{k})_q$ to
$C(\mathfrak{g},\mathfrak{k})_q$.

In the classical case $q=1$ under suitable dominance assumptions on $Z$ the
only non-zero modules $\mathscr{L}_s(Z)$ are those with
$s=\frac12(\dim\mathfrak{k}-\dim(\mathfrak{k}\cap\mathfrak{l}))$ \cite[p.
369]{KV}. So, a particular interest is in considering the
$U_q\mathfrak{g}$-modules
$$
A_\mathfrak{q}(\lambda)\overset{\mathrm{def}}{=}
\mathscr{L}_s(\mathbb{C}_\lambda),\qquad
s=\frac12(\dim\mathfrak{k}-\dim(\mathfrak{k}\cap\mathfrak{l})),
$$
with $\lambda=(\lambda_1,\lambda_2,\ldots,\lambda_l)\in P$ and $\lambda_j=0$
for $j\in\mathbb{L}$.

\renewcommand{\thesection}{A}
\section*{Appendix. Proofs of the Lemmas}

Let $\mathbb{G}\supset\mathbb{L}$ be subsets of $\{1,2,\ldots,l\}$. Equip
$U_q\mathfrak{g}$ with a grading by setting
$$
\deg E_j=
\begin{cases}
1, & j\in\mathbb{G}\setminus\mathbb{L}
\\ 0, & j\in\mathbb{L}
\end{cases},\qquad
\deg F_j=
\begin{cases}
-1, & j\in\mathbb{G}\setminus\mathbb{L}
\\ 0, & j\in\mathbb{L}
\end{cases},
$$
$$\deg K_j^{\pm 1}=0,\qquad j=1,2,\ldots,l.$$

Let $\left(U_q\mathfrak{q}_\mathbb{L}^\pm\right)_j=\left\{\left.\xi\in
U_q\mathfrak{q}_\mathbb{L}^\pm\right|\:\deg\xi=j\right\}$. Obviously,
$\left(U_q\mathfrak{q}_\mathbb{L}^\pm\right)_0=U_q\mathfrak{l}$.

\begin{Lemma}\label{finite_rk}
The homogeneous components $\left(U_q\mathfrak{q}_\mathbb{L}^\pm\right)_j$
are free left and free right $U_q\mathfrak{l}$-modules of finite rank.
\end{Lemma}

{\bf Proof.} Recall that the subsets
$\mathbb{L}\subset\mathbb{G}\subset\{1,2,\ldots,l\}$ determine the
subgroups $W_\mathbb{L}\subset W_\mathbb{G}$ of the Weyl group $W$. Those
subgroups are exactly the Weyl groups of the Lie algebras $\mathfrak{l}$
and $\mathfrak{g}$, respectively. Let $w_{0,\mathbb{L}}$ (respectively,
$w_{0,\mathbb{G}}$) be the longest element of the Weyl group $W_\mathbb{L}$
(respectively, $W_\mathbb{G}$).

We start from the larger subset $\mathbb{G}$. Choose a reduced expression
\begin{equation}\label{re0}
w_{0,\mathbb{G}}=s_{i_1}s_{i_2}s_{i_3}\ldots s_{i_M},
\end{equation}
with $s_{i_k}$ being the reflection corresponding to the simple root
$\alpha_{i_k}$. Consider the Lusztig automorphisms $T_i$, $i\in\mathbb{G}$,
of the algebra $U_q\mathfrak{g}'$ generated by $E_i$, $F_i$, $K_i^{\pm 1}$,
$i\in\mathbb{G}$, \cite[Chapter 8]{Jant}. They have the form
\begin{equation}\label{Lusz}
\begin{gathered}
T_i(K_j)=K_jK_i^{-a_{ij}},
\\ T_i(E_i)=-F_iK_i,\qquad T_i(F_i)=-K_i^{-1}E_i,
\\ T_i(E_j)=\sum_{r+s=-a_{ij}}\operatorname{const}E_i^sE_jE_i^r,\qquad i\ne
j,
\\ T_i(F_j)=\sum_{r+s=-a_{ij}}\operatorname{const}F_i^rF_jF_i^s,\qquad i\ne
j.
\end{gathered}
\end{equation}
In a similar way, consider the algebra $U_q\mathfrak{l}'$ associated to
$\mathbb{L}$. Let $W^\mathbb{L}=\{w\in W_\mathbb{G}|\:l(ws)>l(w)\;\text{for
all simple reflections}\;s\in W_\mathbb{L}\}$, where $l(w)$ stands for the
length of $w$. One has $w_{0,\mathbb{G}}=w^\mathbb{L}w_{0,\mathbb{L}}$ with
$w^\mathbb{L}\in W^\mathbb{L}$ and
$l(w_{0,\mathbb{G}})=l(w^\mathbb{L})+l(w_{0,\mathbb{L}})$ \cite[Proposition
1.10(c)]{Hum2}. This allows one to rewrite the reduced expression
\eqref{re0} in the form
\begin{equation}\label{res}
w_{0,\mathbb{G}}=s_{i_1}s_{i_2}s_{i_3}\ldots
s_{i_{M'}}s_{i_{M'+1}}s_{i_{M'+2}}\ldots s_{i_M},
\end{equation}
with $w_{0,\mathbb{L}}=s_{i_{M'+1}},s_{i_{M'+2}},\ldots,s_{i_M}$ and
$w^\mathbb{L}=s_{i_1}s_{i_2}s_{i_3}\ldots s_{i_{M'}}$.

We are about to apply the Lusztig theorem \cite[Theorem 8.24]{Jant}. Observe
that the special reduced expression \eqref{res}, in view of the explicit
form of the Lusztig automorphisms \eqref{Lusz} implies that all the
monomials
$$
T_{i_1}T_{i_2}\cdots T_{i_{M'-1}}(E_{i_{M'}}^{a_{M'}})\cdot
T_{i_1}T_{i_2}\cdots T_{i_{M'-2}}(E_{i_{M'-1}}^{a_{M'-1}})\cdot\ldots\cdot
T_{i_1}T_{i_2}(E_{i_3}^{a_3})\cdot T_{i_1}(E_{i_2}^{a_2})\cdot E_{i_1}^{a_1}
$$
are in $U_q\mathfrak{l}'$. It follows from the Lusztig theorem \cite[Theorem
8.24]{Jant} that the monomials
\begin{multline}\label{PBW}
T_{i_1}T_{i_2}\cdots T_{i_{M-1}}(E_{i_M}^{a_M})\cdot T_{i_1}T_{i_2}\cdots
T_{i_{M-2}}(E_{i_{M-1}}^{a_{M-1}})\cdot\ldots
\\ \ldots\cdot T_{i_1}T_{i_2}\cdots
T_{i_{M'+1}}(E_{i_{M'+2}}^{a_{M'+2}})\cdot T_{i_1}T_{i_2}\cdots
T_{i_{M'}}(E_{i_{M'+1}}^{a_{M'+1}})
\end{multline}
with all $a_{i_k}\in\mathbb{Z}_+$, form a free basis in the right
$U_q\mathfrak{l}$-module $U_q\mathfrak{q}_\mathbb{L}^+$. As one can readily
separate out for each $j\in\mathbb{Z}_+$ finitely many such monomials that
span the $j$-th homogeneous component, we thus get our claim for
$(U_q\mathfrak{q}_\mathbb{L}^+)_j$ as a right $U_q\mathfrak{l}$-module. All
other claims can be proved in a similar way. \hfill $\square$

\medskip

\begin{Corollary}\label{V_l}
Every $U_q\mathfrak{g}$-module $V$ contains the largest submodule
$V_\mathfrak{l}$ of the category $C(\mathfrak{g},\mathfrak{l})_q$.
\end{Corollary}

{\bf Proof.} Obviously, $V$ possesses the largest weight submodule
$V_\mathfrak{h}$. What remains is to prove that the subspace
$V_\mathfrak{l}=\{v\in V_\mathfrak{h}|\:\dim(U_q\mathfrak{l}v)<\infty\}$ is
a submodule of the $U_q\mathfrak{g}$-module $V$. Let
$\xi\in\left(U_q\mathfrak{q}_\mathbb{L}^\pm\right)_j$, $v\in
V_\mathfrak{l}$. It follows from Lemma \ref{finite_rk} that for some
$\{\eta_1,\eta_2,\ldots,\eta_{N(j)}\}$
$$
\left(U_q\mathfrak{q}_\mathbb{L}^\pm\right)_j=
\sum_{k=1}^{N(j)}\eta_kU_q\mathfrak{l}.
$$
Hence,
$$
\dim(U_q\mathfrak{l}\,\xi
v)\le\dim\left(\left(U_q\mathfrak{q}_\mathbb{L}^\pm\right)_jv\right)=
\dim\left(\sum_{k=1}^{N(j)}\eta_kU_q\mathfrak{l}\,v\right)<\infty.
$$
Thus $\xi v\in V_\mathfrak{l}$. \hfill $\square$

\medskip

\begin{Corollary}\label{g_l}
For any $V\in C(\mathfrak{l},\mathfrak{l})_q$, the module $P(V)$ belongs to
the category $C(\mathfrak{g},\mathfrak{l})_q$.
\end{Corollary}

{\bf Proof.} Since $U_q\mathfrak{g}$ is a free right
$U_q\mathfrak{l}$-module, one has an embedding of $U_q\mathfrak{l}$-modules
$$i:V\hookrightarrow P(V),\qquad i:v\mapsto 1\otimes v.$$
The relation $P(V)=P(V)_\mathfrak{l}$ is due to
$$
P(V)=U_q\mathfrak{g}\,i(V)\subset
U_q\mathfrak{g}\,P(V)_\mathfrak{l}=P(V)_\mathfrak{l}.\eqno\square
$$

\medskip

The rest of the statements of Lemma \ref{P_V} can be proved in the same way
as in the classical case $q=1$.

To prove Lemma \ref{hard_lemma} we need the following auxiliary statement.

\begin{Lemma}\label{easy_lemma}
Let $\lambda\in P_+^\mathbb{L}$ and
$\left(U_q\mathfrak{q}_\mathbb{L}^\pm\right)_j$ be the homogeneous
components of the graded algebras $U_q\mathfrak{q}_\mathbb{L}^\pm$.

1. The vector spaces
$\left(U_q\mathfrak{q}_\mathbb{L}^\pm\right)_j\cdot\mathscr{P}_\lambda
\subset\operatorname{End}V^\mathrm{univ}$ are finite dimensional.

2. The vector spaces
$\mathscr{P}_\lambda\cdot\left(U_q\mathfrak{q}_\mathbb{L}^\pm\right)_j
\subset\operatorname{End}V^\mathrm{univ}$ are finite
dimensional.\footnote{A dot is used here to denote the product of elements
in $\operatorname{End}V^\mathrm{univ}$.}
\end{Lemma}

{\bf Proof.} Prove the first statement. We consider
$\operatorname{End}V^\mathrm{univ}$ as a $U_q\mathfrak{l}$-module with
respect to the action as follows:
$$
(\xi a):v\mapsto\xi(av),\qquad v\in V^\mathrm{univ},\quad
a\in\operatorname{End}V^\mathrm{univ},\quad\xi\in U_q\mathfrak{l}.
$$
Obviously,
$\left(U_q\mathfrak{q}_\mathbb{L}^\pm\right)_j\cdot\mathscr{P}_\lambda$ is a
submodule of the $U_q\mathfrak{l}$-module
$\operatorname{End}V^\mathrm{univ}$.

Let
$\pi_\lambda:U_q\mathfrak{l}\to\operatorname{End}L(\mathfrak{l},\lambda)$ be
the representation of the Hopf subalgebra $U_q\mathfrak{l}$ corresponding to
the $U_q\mathfrak{l}$-module $L(\mathfrak{l},\lambda)$. If
$\xi\in\operatorname{Ker}\pi_\lambda$, then $\xi\cdot\mathscr{P}_\lambda=0$.
Hence the diagram
$$
\xymatrix{U_q\mathfrak{l}\ar[d]_{\pi_\lambda}
\ar[r]^-{\xi\mapsto\xi\cdot\mathscr{P}_\lambda} &
\operatorname{End}V^\mathrm{univ}\\
\operatorname{End}L(\mathfrak{l},\lambda)\ar[ru]}
$$
can be completed up to a commutative one, and
$$
\dim(U_q\mathfrak{l}\cdot\mathscr{P}_\lambda)\le
(\dim(L(\mathfrak{l},\lambda)))^2.
$$
Thus the first statement of the Lemma is proved in the case $j=0$. It
remains to elaborate the fact that $(U_q\mathfrak{q}_\mathbb{L}^\pm)_j$ is a
right $U_q\mathfrak{l}$-module of finite rank (see Lemma \ref{finite_rk}).

Now the first statement is proved. The second one can be proved in a similar
way using the commutative diagram
$$
\xymatrix{U_q\mathfrak{l}\ar[d]_{\pi_\lambda}
\ar[r]^-{\xi\mapsto\mathscr{P}_\lambda\cdot\xi} &
\operatorname{End}V^\mathrm{univ}\\
\operatorname{End}L(\mathfrak{l},\lambda)\ar[ru]}
$$
in the category of $U_q\mathfrak{l}^\mathrm{op}$-modules.\footnote{The
algebra $U_q\mathfrak{l}^\mathrm{op}$ is derived from $U_q\mathfrak{l}$ by
replacing the multiplication with the opposite one.}\hfill $\square$

\medskip

Turn to the proof of Lemma \ref{hard_lemma}. Consider the subspace $A\subset
U_q\mathfrak{g}$ of such elements $\xi\in U_q\mathfrak{g}$ that for any
$\lambda\in P_+^\mathbb{L}$ there exists a finite subset $\Lambda\subset
P_+^\mathbb{L}$ with
\begin{eqnarray}
\xi\cdot\mathscr{P}_\lambda &=&
\chi_\Lambda\cdot\xi\cdot\mathscr{P}_\lambda,\label{fp}
\\ \mathscr{P}_\lambda\cdot\xi &=&
\mathscr{P}_\lambda\cdot\xi\cdot\chi_\Lambda.\label{ip}
\end{eqnarray}
One has to prove that $A=U_q\mathfrak{g}$. For that, it suffices to
demonstrate that $A$ is a subalgebra of $U_q\mathfrak{g}$ and
$A\supset\left(U_q\mathfrak{q}_\mathbb{L}^\pm\right)_j$ for all $j$.

Prove the first statement. It is obvious that $A$ is a vector subspace. Now
let $\xi,\eta\in A$, and consider the product $\xi\eta$. A double
application of \eqref{fp} first w.r.t. $\eta$ and then $\xi$ allows one to
deduce that the image of the linear operator
$\xi\cdot\eta\cdot\mathscr{P}_\lambda\in\operatorname{End}V^\mathrm{univ}$
is accommodated by the sum of finitely many isotypic components
$V_\lambda^\mathrm{univ}$. Let this sum be
$\bigoplus\limits_{\lambda'\in\Lambda'}V_{\lambda'}^\mathrm{univ}$, then
$\xi\cdot\eta\cdot\mathscr{P}_\lambda=
\chi_{\Lambda'}\cdot\xi\cdot\eta\cdot\mathscr{P}_\lambda$. Thus we get
\eqref{fp} for $\xi\eta$.

In a similar way, apply \eqref{ip} twice to deduce that for some finite
subset $\Lambda''\subset P_+^\mathbb{L}$ one has
$\operatorname{Ker}(\mathscr{P}_\lambda\cdot\xi\cdot\eta)\supset
\left(\mathrm{id}-\chi_{\Lambda''}\right)V^\mathrm{univ}$, hence
$\mathscr{P}_\lambda\cdot\xi\cdot\eta=\mathscr{P}_\lambda\cdot\xi\cdot\eta
\cdot\chi_{\Lambda''}$. Thus \eqref{ip} holds for $\xi\eta$.

Turn to the second statement. We restrict ourselves to proving that
\eqref{fp} is valid for all
$\xi\in\left(U_q\mathfrak{q}_\mathbb{L}^\pm\right)_j$. It follows from
Lemma \ref{easy_lemma} that for every $j$, $\lambda$, the vector space
$\left(U_q\mathfrak{q}_\mathbb{L}^\pm\right)_j\cdot\mathscr{P}_\lambda$ is
finite dimensional. It is also a $U_q\mathfrak{l}$-module since
$U_q\mathfrak{l}\left(U_q\mathfrak{q}_\mathbb{L}^\pm\right)_j\subset
\left(U_q\mathfrak{q}_\mathbb{L}^\pm\right)_j$. Hence
$\left(U_q\mathfrak{q}_\mathbb{L}^\pm\right)_j\cdot\mathscr{P}_\lambda$ is
a sum of finitely many $U_q\mathfrak{l}$-isotypic components. On the other
hand, if for some $a\in\operatorname{End}V^\mathrm{univ}$, the
$U_q\mathfrak{l}$-module $U_q\mathfrak{l}\cdot
a\subset\operatorname{End}V^\mathrm{univ}$ is a multiple of
$L(\mathfrak{l},\mu)$, $\mu\in P_+^\mathbb{L}$, one has
$aV^\mathrm{univ}\subset V_\mu^\mathrm{univ}$ and $\mathscr{P}_\mu\cdot
a=a$. It follows that
$\mathscr{P}_{\lambda'}\cdot\left(U_q\mathfrak{q}_\mathbb{L}^\pm\right)_j
\cdot\mathscr{P}_\lambda=0$ for all but finitely many $\lambda'\in
P_+^\mathbb{L}$. Hence
$\xi\cdot\mathscr{P}_\lambda=\chi_\Lambda\cdot\xi\cdot\mathscr{P}_\lambda$
for some finite subset $\Lambda\subset P_+^\mathbb{L}$. \hfill $\square$

\end{document}